\documentclass[10pt,a4paper]{elsarticle}

\usepackage{color}
\usepackage[utf8]{inputenc}
\usepackage[english]{babel}
\usepackage{hyperref}
\usepackage{amsmath,amsthm,enumerate}
\usepackage{amssymb}
\usepackage{graphicx}
\usepackage{tikz}
\usetikzlibrary{calc,decorations.pathmorphing,decorations.text,fixedpointarithmetic}
\usepackage{subfigure}
\usepackage{refcount}
\usepackage{cleveref}
\usepackage{MnSymbol}
\usepackage[hmargin=2.5cm,vmargin=3cm]{geometry}

\newtheorem{theorem}{Theorem}
\newtheorem{proposition}[theorem]{Proposition}
\newtheorem{lemma}[theorem]{Lemma}
\newtheorem{observation}[theorem]{Observation}
\newtheorem{corollary}[theorem]{Corollary}

\newtheorem{conjecture}[theorem]{Conjecture}
\newtheorem{example}[theorem]{Example}

\newtheorem{problem}[theorem]{Problem}

\theoremstyle{definition}
\newtheorem{definition}[theorem]{Definition}

\begin{document}
	
	\begin{frontmatter}
		\title{Bounding signed bipartite partial $t$-trees and application to edge-coloring}

		\author[XMUT]{Meirun Chen}
		\author[IRIF,ZNU]{Reza Naserasr}
		
		\address[XMUT]{School of Mathematics and Statistics, Xiamen University of Technology, Xiamen Fujian 361024, China. mrchen@xmut.edu.cn.}
		\address[IRIF]{Université de Paris, IRIF, CNRS, F-75013 Paris, France. reza@irif.fr}
	    \address[ZNU]{School of Mathematical Sciences, Zhejiang Normal University, Jinhua {\rm 321004}, China}

		\begin{abstract}
			Given a signed bipartite graph $(B, \pi)$ of negative girth $2k$, we present a necessary and sufficient condition for it to have the following property: each signed bipartite graph $(G, \sigma)$ whose negative girth is at least $2k$ and whose underlying graph has treewidth at most $t$ admits a homomorphism to $(B, \pi)$. 
			
			Applying the result on the signed projective cube $SPC(2k-1)$, we conclude that every signed bipartite graph of negative girth at least $2k$ whose underlying graph is a partial 3-tree admits a homomorphism to $SPC(2k-1)$. For planar partial 3-trees, applying duality we conclude that if $G$ is a planar $2k$-regular multigraph whose dual has treewidth at most 3 and such that every edge-cut $(X, V\backslash X)$, where $|X|$ is odd, has size at least $2k$, then $G$ is $2k$-edge-colorable. This supports a conjecture of Seymour which, in full generality, largely extends Tait's reformulation of the four-color theorem, claiming that the fractional edge-chromatic number of a planar multigraph determines its edge-chromatic number.
			
			Finally, noting the contrast between fractional isomorphism and quantum isomorphism, where the former admits a polynomial time algorithm while the latter is proved to be undecidable, and observing the similarities of these notions to the subject of our study, we ask if there is an algorithm to decide if an input signed graph $\widehat{B}$ has the following property: if a signed planar graph $\widehat{G}$ does not map to $\widehat{B}$, it would be because a cycle in $\widehat{G}$ does not map to $\widehat{B}$. In other words, minimal planar graphs that do not map to $\widehat{B}$ are signed cycles.
 	\end{abstract}
		
		\begin{keyword}
			Signed graphs, homomorphism, partial tree, edge-coloring
		\end{keyword}
		
	\end{frontmatter}

	\section{Introduction}
	
	The chromatic number of a (loop-free) graph can be, equivalently, defined as the minimum number of vertices in a loop-free homomorphic image. Motivated by the question of 4-coloring maps, bounding the chromatic number of classes of graphs then is one of the most central areas of graph theory. One of the earliest observations in the study of graph homomorphisms is that: there is always an odd cycle in the homomorphic image, $f(C)$, of any odd cycle $C$ whose length is (naturally) at most the length of $C$. In particular, the image of every loop (a 1-cycle) is a loop. 
	
	Generalizing the notion of the chromatic number from this point of view, one may ask: Given a graph $G$ of odd-girth at least $2k+1$ what is the smallest homomorphic image of $G$ whose odd-girth is also $2k+1$? Many of the fundamental questions in graph theory, including the notion of the chromatic number, can be viewed as special cases of this question.
	
	Furthering this question, given a graph $B$ of odd-girth $2k+1$, and a class $\mathcal{C}$ of graphs of odd-girth at least $2k+1$, one may ask whether every graph in $\mathcal{C}$ maps to $B$. When $\mathcal{C}$ is the class of graphs of treewidth at most $t$, a necessary and sufficient condition for $B$ to be an affirmative case of this question is answered in \cite{CN20}.
	
	To better understand the correlation between homomorphisms and minor theory, the notion of homomorphisms of signed graphs has been developed, see \cite{NRS15} and \cite{NSZ21}. In this work then we address the analogous question for signed bipartite graphs. More precisely we ask:
	
	\begin{problem} 
		Given a signed bipartite graph $(B,\pi)$ of negative girth $2k$, does every signed bipartite partial $t$-tree of negative girth at least $2k$ admit a homomorphism to $(B, \pi)$. What is the smallest possible order of such a $(B, \pi)$ in terms of $k$ and $t$? 
	\end{problem}

	In this work, we provide a necessary and sufficient condition for $(B,\pi)$ to receive an affirmative answer. Furthermore, when $t=3$ we show that signed projective cube of dimension $2k-1$, a signed graph on $2^{2k-1}$ vertices, is an optimal solution. 
	To better address the question, we shall first settle the notation and introduce the necessary concepts.

	\section{Preliminaries}
	
	\subsection{Signed graphs}
	
	A signed graph $(G,\sigma)$ is a graph $G$ together with a signature function $\sigma:E(G)\rightarrow\{+,-\}$. A \emph{switching} at a vertex $v$ of a signed graph $(G, \sigma)$ is to multiply the sign of each edge incident with $v$ by $-$, noting that a loop at $v$ is considered to have both its ends at $v$ and will not change sign then. When signature $\sigma$ is of no importance, we may write $\widehat{G}$ instead of $(G, \sigma)$.
	
	The sign of a structure $S$ of a signed graph $(G,\sigma)$ is the product of the signs of all the  edges in $S$, considering multiplicities. For some structures such as a cycle, or a closed walk in general, a switching will not affect the sign. For other structures, a switching may or may not change the sign. For example, for a path $P$ a switching at one end-vertex will change the sign of $P$ while a switching at an internal-vertex of $P$ will not affect its sign. 
	
	For a pair $x,y$ of vertices of a signed graph $(G, \sigma)$ their \emph{distance}, denoted $d_G(x,y)$, or simply $d(x,y)$ when $G$ is clear from the context, is their distance in $G$. Their \emph{algebraic distance}, denoted $ad_{(G, \sigma)}(x,y)$ or simply $ad(x,y)$, is $d(x,y)\sigma(P)$, where $P$ is a shortest $x-y$ path. When there are both positive and negative $x-y$ paths of the shortest length, then we choose a positive one. Observe that a switching at $x$ changes the sign of almost all algebraic distances to $x$. To be precise, when there are both positive and negative paths of length $d(x,y)$, then the sign does not change. However, this case is just a matter of settling notation and will not impact the study. 
	
	A \emph{negative cycle} (respectively, \emph{positive cycle}) is commonly referred to as \emph{unbalanced}  (respectively, \emph{balanced}) cycle. A negative cycle of length $l$ is denoted by $C_{-l}$. The length of a shortest negative cycle of $(G, \sigma)$ is called the \emph{unbalanced girth} or \emph{negative girth} of $(G, \sigma)$ and is denoted by $g_{-}(G, \sigma)$.

	A fundamental property of signed graphs, which distinguishes it from 2-edge-colored graphs, is the following result of Zaslavsky \cite{Z82}.
	
	\begin{theorem}
		A signed graph $(G,\sigma_1)$ can be obtained from $(G,\sigma_2)$ by switching operations if and only if $\sigma_1(C)=\sigma_2(C)$ for any cycle $C$ of $G$.
	\end{theorem}
	
	\subsection{Homomorphisms}
	
	A \emph{homomorphism} of a graph $G$ to a graph $H$ is a mapping $f:V(G)\rightarrow V(H)$ such that if $xy\in E(G)$ then $f(x)f(y)\in E(H)$. 
	A \emph{homomorphism} of a signed graph $(G,\sigma)$ to a signed graph $(H,\pi)$ is a homomorphism $f$ of $G$ to $H$ such that $\sigma(C)=\pi(f(C))$ for any cycle $C$ in $G$.
	That is equivalent to a switching $(G,\sigma')$ of $(G,\sigma)$ and a mapping of $(G,\sigma')$ to $(H,\pi)$ which not only preserves adjacency but also preserves the signs of the edges.

	\begin{definition} We say a signed graph $\widehat{B}$ \emph{bounds} a class of signed graphs $\mathcal{SG}$,  if there exists a homomorphism from each member of $\mathcal{SG}$ to $\widehat{B}$.
	\end{definition}

	Given a signed graph $(G, \sigma)$ and a pair $x,y$ of its vertices we define the positive distance of $x$ and $y$, denoted $d^+(x,y)$, to be the length of a shortest ($x-y$)-path whose sign is positive. The negative distance, denoted $d^-(x,y)$, is defined analogously. We may observe that $d^+(x,y)$ and $d^-(x,y)$ are not invariant under switching at $x$ (or $y$), but such a switching will interchange the two values. To better distinguish the values, $d^-(x,y)$ will be presented by negative number. This view leads to extension of the notion of signed graphs to signed edge-weighted graphs which is an essential part of this work, see Section 3.  
	
	\subsection{Partial $t$-trees}
	
	A \emph{$t$-tree} is any graph in the class $\mathcal{T}_t$ of graphs defined as follows:
	
	\begin{itemize}
		\item $K_{t+1}$ is in  $\mathcal{T}_t$.
		\item If $G$ is in  $\mathcal{T}_t$ and $K$ is a subgraph isomorphic to $K_t$, then the graph obtained from $G$ by adding a vertex which is joined to all vertices of $K$ is in  $\mathcal{T}_t$. 
	\end{itemize} 
	
	It is easily observed that  $\mathcal{T}_1$ is the class of all trees, hence the name $t$-trees. Any subgraph of a graph $G$ in  $\mathcal{T}_t$ is called a \emph{partial $t$-tree}, the class of which is denoted $\mathcal{PT}_t$. It is a classic homework that a graph $G$ is a partial $t$-tree if and only if its treewidth is at most $t$. Thus we may take this as the definition of treewidth as well.
	
	A $t$-tree is thus built by starting with a $(t+1)$-clique on vertices $v_1, \ldots, v_{t+1}$, which we denote by $H_{t+1}$. Then when adding the vertex $v_{j}$, $j\geq t+2$, we create a new $(t+1)$-clique which we denote by $H_{j}$. The sequence $H_{t+1}, H_{t+2},\ldots,H_{t+l}$  is referred to as \emph{$(t+1)$-cliques sequence} of $G$. Observe that a $t$-tree $G$ can be presented as a $t$-tree with different orderings of vertices. Each ordering leads to a different sequence of $(t+1)$-cliques, but the set of these cliques remain the same and is independent of the orderings.
	
	Given a partial $t$-tree $G$, a $(t+1)$-cliques sequence of $G$ is a $(t+1)$-cliques sequence of a $t$-tree $G'$ which contains $G$ as a spanning subgraph. That is to say, we add as many edges as possible to $G$ so that it becomes a $t$-tree and then take a $(t+1)$-clique sequence of it.
	
	Since one may add different sets of edges in order to get a $t$-tree from a partial $t$-tree, even the set of $(t+1)$-cliques of a partial $t$-tree is not unique.  
	
\begin{figure}[ht]
		\centering 
		
		\begin{minipage}{0.24\textwidth}
			\centering
			\begin{tikzpicture}
				\foreach \i in {1,...,5}
				{
					\draw[rotate=-72*(\i)] (0, 1.5) node[circle, fill=white, draw=black!80, inner sep=.7mm, minimum size=1mm]  (\i){};
				}
				
				\foreach \i/\j in {1/2, 1/3, 1/4, 1/5, 2/3, 2/4, 2/5, 3/4, 3/5, 4/5, 5/1}
				{
					
					\draw[thick] (\i) -- (\j);
				}

			\end{tikzpicture}
			
		\end{minipage}
		\begin{minipage}{0.24\textwidth}
			\centering
			\begin{tikzpicture}
				\foreach \i in {1,2,3}
				{
					\draw[rotate=-120*(\i)] (0, 0.9) node[circle, fill=white, draw=black!80, inner sep=0.7mm, minimum size=1mm]  (\i){};
				}
				
				\foreach \i/\j in {1/2, 1/3, 2/3}
				{
					
					\draw[thick] (\i) -- (\j);
				}
				
				\foreach \i in {4,5,6}
				{
					\draw[rotate=-120*(\i)] (0, 1.8) node[circle, fill=white, draw=black!80, inner sep=0.7mm, minimum size=1mm]  (\i){};
				}
				
				\foreach \i/\j in {4/5, 4/6, 5/6, 1/5, 2/6, 3/4, 1/4, 2/5, 3/6}
				{
					
					\draw[thick] (\i) -- (\j);
				}

			\end{tikzpicture}
			
		\end{minipage}
		\begin{minipage}{0.24\textwidth}
			\centering
			\begin{tikzpicture}
				\foreach \i in {1,...,8}
				{
					\draw[rotate=-45*(\i)] (0, 1.5) node[circle, fill=white, draw=black!80, inner sep=0.7mm, minimum size=1mm]  (\i){};
				}
				
				\foreach \i/\j in {1/2, 2/3, 3/4, 4/5, 5/6, 6/7, 7/8, 8/1, 1/5, 2/6, 3/7, 4/8}
				{
					
					\draw[thick] (\i) -- (\j);
				}

			\end{tikzpicture}
			
		\end{minipage}
		\begin{minipage}{0.24\textwidth}
			\centering
			\begin{tikzpicture}
				\foreach \i in {1,...,5}
				{
					\draw[rotate=-72*(\i)] (0, 1.5) node[circle, fill=white, draw=black!80, inner sep=0.7mm, minimum size=1mm]  (\i){};
				}
				
				\foreach \i/\j in {1/2, 2/3, 3/4, 4/5, 5/1}
				{
					
					\draw[thick] (\i) -- (\j);
				}
				
				\foreach \i in {6,...,10}
				{
					\draw[rotate=-72*(\i)] (0, 0.8) node[circle, fill=white, draw=black!80, inner sep=0.7mm, minimum size=1mm]  (\i){};
				}
				
				\foreach \i/\j in {6/7, 7/8, 8/9, 9/10, 6/10, 1/6, 2/7, 3/8, 4/9, 5/10}
				{
					
					\draw[thick] (\i) -- (\j);
				}

			\end{tikzpicture}
			
		\end{minipage}
		\caption{Forbidden minors of the class of partial 3-trees.}
		\label{fig:forbiddenminorsofpartial3trees}
		
	\end{figure}
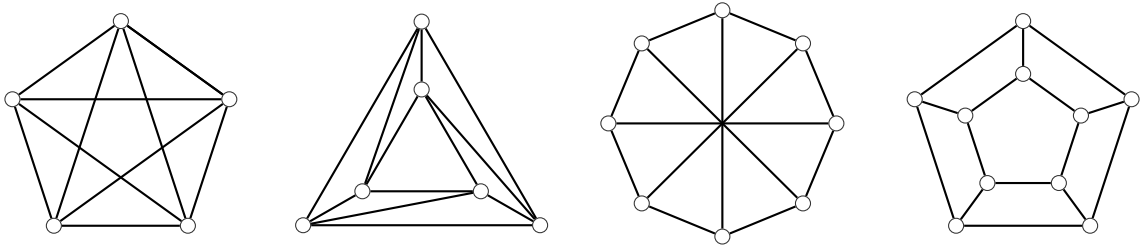

    It is a well known fact that $\mathcal{PT}_t$ is a minor closed family of graphs. In particular $\mathcal{PT}_1$, the class of forests, is the class of graphs with no $K_3$-minor. The class $\mathcal{PT}_2$, is the same as the class of series parallel graphs which is also the same as graphs with no $K_4$-minor. The class $\mathcal{PT}_3$ is the class of graphs with four graphs of \Cref{fig:forbiddenminorsofpartial3trees} as forbidden minors. Thus $\mathcal{PT}_3$ is a subclass of $K_5$-minor free graphs, but contains elements that are not planar.

	The class of outerplanar graphs is a subclass of $\mathcal{PT}_2$ which itself is included in $\mathcal{PT}_3$. These three classes of graphs are natural subclass of $K_5$-minor free graphs which are commonly used to test conjectures generalizing the four-color theorem. The class $\mathcal{PT}_3$, which is the largest of the three, is of special interest in this work.

	\section{Signed edge-weighted graphs}

	A \emph{signed edge-weighted graph} $(H,w)$ is a graph equipped with a mapping $w: E(H)\rightarrow \mathbb{Z}-\{0\}$ which maps each edge $e$ of $H$ to a nonzero integer $w(e)$.

	A signed edge-weighted graph $(H,w)$ is bipartite if there exists a partition $\{X,Y\}$ of $V(H)$ such that for any edge $e$ with one endpoint in $X$ and the other endpoint in $Y$, $w(e)$ is odd and for any edge $e$ with both endpoints in $X$ (resp. $Y$), $w(e)$ is even. In other words, though there might be odd cycles in the underlying graph $H$, the sum of the weights of edges of each cycle is even.

	The \emph{length} of a cycle $C$ in a signed edge-weighted graph $(H,w)$ is the sum of the absolute values of the weights of the edges of $C$. A cycle is an \emph{unbalanced cycle} if the number of edges with negative weights is odd. The unbalanced-girth of a signed edge-weighted graph $(H,w)$ is the length of a shortest unbalanced cycle.
	
	Given a signed edge-weighted graph $(H,w)$ where $1\leq |w(e)|\leq k$ for every edge $e$, we build a signed graph $\overline{\overline{(H,w)}}_{-2k}$ as follows: each edge $e=xy$ is replaced by two paths connecting $x$ and $y$ whose internal vertices are distinct, one of the paths is of length $|w(e)|$ and the other is of length $2k-|w(e)|$; the assignment of signs to the edges of the two paths is such that the path associated with $w(e)$ has the same sign as $w(e)$ and the other one has the opposite sign. Thus the two paths together form a negative $2k$-cycle. Observe that a switching at a vertex of $(H,w)$ corresponds to the switching of the same vertex in $\overline{\overline{(H,w)}}_{-2k}$. An example of weighted signed $(K_4,w)$ is presented in \Cref{fig:WeightedK4} where the red weights have the negative sign. The corresponding signed graph $\overline{\overline{(K_4,w)}}_{-4}$ is presented in \Cref{fig:barbarGraph}.	 
	It is worth to mention that if the signed edge-weighted graph $(H,w)$ is bipartite, then every cycle of $\overline{\overline{(H,w)}}_{-2k}$ is of even length, that is to say its underlying graph is bipartite.

	\begin{figure}[ht]
		\centering 
		
		\begin{minipage}{0.45\textwidth}
			\centering
			\begin{tikzpicture}
				\foreach \i in {1,2}
				{
					\draw[rotate=-120*(\i)] (0, 3) node[circle, fill=white, draw=black!80, inner sep=.7mm, minimum size=1mm]  (\i){};
					
				}
					\foreach \i in {3}
					{
						\draw[rotate=-120*(\i)] (0, 3) node[circle, fill=lightgray, draw=black!80, inner sep=.7mm, minimum size=1mm]  (\i){};
					}
					
					\draw (0, 0) node[circle, fill=lightgray, draw=black!80, inner sep=.7mm, minimum size=1mm]  (4){};

					\foreach \i in {5}
				{
					\draw[rotate=-120*(\i)] (0.1, 1) node[circle, fill=white, draw=black!80, inner sep=.0mm, minimum size=0mm]  (\i){}; 
					
				}
				
					\foreach \i in {6}
				{
					\draw[rotate=-120*(\i)] (0.1, 1) node[circle, fill=white, draw=black!80, inner sep=.0mm, minimum size=0mm]  (\i){};
					
				}
				
					\foreach \i in {7}
				{
					\draw[rotate=-120*(\i)] (0.1, 1) node[circle, fill=white, draw=black!80, inner sep=.0mm, minimum size=0mm]  (\i){};
					
				}
				
						\foreach \i in {8}
				{
					\draw[rotate=-120*(\i)] (0.1, 1.6) node[circle, fill=white, draw=black!80, inner sep=.0mm, minimum size=0mm]  (\i){};
					
				}

						\foreach \i in {9}
				{
					\draw[rotate=-120*(\i)] (0.1, 1.4) node[circle, fill=white, draw=black!80, inner sep=.0mm, minimum size=0mm]  (\i){};
					
				}
				
						\foreach \i in {10}
				{
					\draw[rotate=-120*(\i)] (0.1, 1.6) node[circle, fill=white, draw=black!80, inner sep=.0mm, minimum size=0mm]  (\i){};
					
				}

				\foreach \i/\j in {1/2, 1/3, 2/3, 1/4, 2/4, 3/4}
				{
					
					\draw[thick] (\i) -- (\j);
				}

					\foreach \i/\j in {9/6, 8/5}
				{
					
					\draw[dashed, thick, blue] (\i) -- (\j);
				}

					\foreach \i/\j in {7/10}
				{
					
					\draw[thick, dashed, red] (\i) -- (\j);
				}

							\foreach \i in {1,3}
				{
					\draw[rotate=-120*(\i)+51] (0, 1.7) node[circle, fill=white, draw=black!80, inner sep=.0mm, minimum size=0mm]  (\i){}; 
					
				}

							\foreach \i in {5,6}
				{
					\draw[rotate=-120*(\i)-51] (0, 1.7) node[circle, fill=white, draw=black!80, inner sep=.0mm, minimum size=0mm]  (\i){}; 
					
				}

							\foreach \i in {2}
				{
					\draw[rotate=-120*(\i)+53] (0, 1.7) node[circle, fill=white, draw=black!80, inner sep=.0mm, minimum size=0mm]  (\i){}; 
					
				}

				\foreach \i in {7}
				{
					\draw[rotate=-120*(\i)-53] (0, 1.7) node[circle, fill=white, draw=black!80, inner sep=.0mm, minimum size=0mm]  (\i){}; 
					
				}
				
					\foreach \i/\j in {1/6, 2/7}
				{
					
					\draw[thick, dashed, red] (\i) -- (\j);
				}
				
					\foreach \i/\j in {3/5}
				{
					
					\draw[thick, dashed, blue] (\i) -- (\j);
				}
				
					\draw (0, -2) node[circle, fill=white, draw=black!80, inner sep=0mm, minimum size=0mm]  (){};
			\end{tikzpicture}
				\caption{A bipartite $(K_4,w)$, red weights are negative.}
			\label{fig:WeightedK4}
		\end{minipage}
		\begin{minipage}{0.45\textwidth}
			\centering
			\begin{tikzpicture}
				\foreach \i in {1,2,3}
				{
					\draw[rotate=-120*(\i)] (0, 3) node[circle, fill=white, draw=black!80, inner sep=0.7mm, minimum size=1mm]  (\i){};
					\draw (0, 0) node[circle, fill=white, draw=black!80, inner sep=0.7mm, minimum size=1mm]  (4){};
					\draw (0, 1.5) node[circle, fill=white, draw=black!80, inner sep=.5mm, minimum size=1mm]  (5){};
					\draw (0, -1.5) node[circle, fill=white, draw=black!80, inner sep=.5mm, minimum size=1mm]  (6){};
					\draw (0.84, -.5) node[circle, fill=white, draw=black!80, inner sep=.5mm, minimum size=1mm]  (7){};
					\draw (1.7, -1) node[circle, fill=white, draw=black!80, inner sep=.5mm, minimum size=1mm]  (8){};
					\draw (-0.84, -.5) node[circle, fill=white, draw=black!80, inner sep=.5mm, minimum size=1mm]  (9){};
					\draw (-1.7, -1) node[circle, fill=white, draw=black!80, inner sep=.5mm, minimum size=1mm]  (10){};
					\draw (0, -2) node[circle, fill=white, draw=black!80, inner sep=.5mm, minimum size=1mm]  (11){};
					\draw (0.86, 1.5) node[circle, fill=white, draw=black!80, inner sep=.5mm, minimum size=1mm]  (12){};
					\draw (1.74, 0) node[circle, fill=white, draw=black!80, inner sep=.5mm, minimum size=1mm]  (13){};
					\draw (-0.86, 1.5) node[circle, fill=white, draw=black!80, inner sep=.5mm, minimum size=1mm]  (14){};
					\draw (-1.74, 0) node[circle, fill=white, draw=black!80, inner sep=.5mm, minimum size=1mm]  (15){};
					\draw (-0.3, 1.5) node[circle, fill=white, draw=black!80, inner sep=.5mm, minimum size=1mm]  (16){};
				}
				
				\foreach \i/\j in {3/5, 4/5, 1/6, 8/1, 4/7, 4/9, 9/10, 2/10, 12/13, 3/12, 3/14, 2/15, 14/15}
				{
					
					\draw[thick, blue] (\i) -- (\j);
				}
				
				\foreach \i/\j in {8/7, 1/13, 2/6}
				{
					
					\draw[thick, red] (\i) -- (\j);
				}
				
				\draw[bend left=20, thick, red] (4) to (2);
				\draw[bend right=20, thick, red] (3) to (2);
				\draw[bend left=20, thick, blue] (4) to (1);
				\draw[bend left=20, thick, blue] (3) to (1);
				\draw[bend right=10, thick, blue] (11) to (1);
				\draw[bend left=10, thick, blue] (11) to (2);
				\draw[bend left=10, thick, red] (16) to (3);
				\draw[bend right=10, thick, blue] (16) to (4);
				
			\end{tikzpicture}
			
				\caption{$\overline{\overline{(K_4,w)}}_{-4}$.}
			\label{fig:barbarGraph}
		\end{minipage}
	
	\end{figure}
	
		In the following definition, we extend the notion of negative girth to signed edge-weighted graphs in some sense.
	
	\begin{definition}\label{lem:g-wide} 
		A signed edge-weighted graph $(H,w)$, where $1\leq |w(e)|\leq k$ for every edge $e$, is said to be $2k$-wide if the signed graph $\overline{\overline{(H,w)}}_{-2k}$ is of unbalanced girth $2k$. 
	\end{definition}
	
	\begin{example}\label{exm:all-k-weight}
		Let $w$ be an assignment of weights to the edges of $K_r$ such that $w(e)\in\{k,k-1\}$ and the sum of the weights of each triangle is even. Then $(K_r,w)$ is bipartite and $2k$-wide.
	\end{example}
	
	In this example, the fact that weights of all triangles are even implies that the (signed) edge-weighted graph itself is bipartite. Assuming that $X\cup Y$ is the bipartition of vertices, it follows that the weight of each edge with both ends in the same part is the even value among $k-1$ and $k$ and the weight of each edge connecting the two parts is the odd value among the two.  
	
	Necessary and sufficient condition for a  bipartite signed edge-weighted triangle $(K_{3},w)$ to be $2k$-wide was given in  \cite{BFN19}.
	
	\begin{proposition}\label{prop:2k+1WideTriangleTest1}
		Let $(K_{3},w)$ be a bipartite signed edge-weighted triangle with edge weights $a,b,c$ satisfying: $a+b+c$ is even, and $1\leq|a|\leq |b|\leq |c|\leq k$. Then $(K_{3},w)$ is $2k$-wide if and only if either:\\
		$(i)$ $abc<0$ and $|a|+|b|+|c|\geq 2k$, or\\
		$(ii)$ $abc>0$ and $|a|+|b|\geq |c|$.
	\end{proposition}

	Next we provide a necessary and sufficient condition for a bipartite signed edge-weighted complete graph $(K_{t},w)$ to be $2k$-wide. 

	\begin{theorem}\label{thm:complete weighted graph} A bipartite signed edge-weighted complete graph $(K_{t},w)$ satisfying $|w(e)|\leq k$ for every edge $e$, is $2k$-wide if and only if each of its triangles is $2k$-wide. 
	\end{theorem}
	
	\begin{proof} The necessary part is trivial; if $(K_{t},w)$ is $2k$-wide, then any of its subgraphs is also $2k$-wide.
		
		Suppose that each triangle of the bipartite signed edge-weighted complete graph $(K_t,w)$ is $2k$-wide. Now we show that $(K_t,w)$ is also $2k$-wide. That is to say, the signed graph $\overline{\overline{(K_t,w)}}_{-2k}$ is of unbalanced girth $2k$. Towards a contradiction assume there are negative cycles of length smaller than $2k$ in $\overline{\overline{(K_t,w)}}_{-2k}$. Among all such cycles, let $C$ be one which has minimum number of vertices from $K_t$.
		
		We first observe that $C$ has at least four vertices of $K_t$. Otherwise, for some triangle $T$ of $K_t$ the cycle $C$ is a subgraph of $\overline{\overline{(T,w)}}_{-2k}$. This contradicts the assumption that every triangle, in particular $T$, is $2k$-wide.

		Let	$u,v,w,x$ be four of the vertices of $K_t$ which are also vertices of $C$, appearing in this cyclic order on $C$. Let $C_{uv}$, $C_{vw}$, $C_{wx}$, and $C_{xu}$ be the segments of the cycle connecting the corresponding vertices and not containing the other two vertices. Let $a$, $b$, $c$, and $d$ be the algebraic length of these segments, respectively. Since $C$ is a negative cycle of length less than $2k$, we have $abcd<0$ and $|a|+|b|+|c|+|d|<2k$. Since $abcd$ is negative, we either have $ab<0$ and $cd>0$ or $ab>0$ and $cd<0$. By symmetries, we assume the former holds. Let $P$ and $Q$ be the two threads connecting $u,w$ in $\overline{\overline{(K_t,w)}}_{-2k}$, noting that one of $P,Q$, say $P$, is positive, the other, $Q$, is negative. Next we consider the following two negative cycles: $PC_{uv}C_{vw}$ and  $QC_{wx}C_{xu}$. The total length of these two (negative) cycles is less than $4k$. Thus one of them has length less than $2k$, but each of them has less vertices from $K_t$, contradicting the choice of $C$. 
	\end{proof}

	\begin{lemma}\label{thm:a new vertex} Let $(H,w)$ be a $2k$-wide bipartite signed edge-weighted  $t$-tree. If we add a new vertex $u$ and connect it to $t$ vertices of a $t$-clique $K_t$ of $H$ and assign a nonzero weight to each new edge such that the resulting $(t+1)$-clique, $(K, w)$, is bipartite and $2k$-wide, then the resulting bipartite signed  edge-weighted $t$-tree is also $2k$-wide.
	\end{lemma}
	
	\begin{proof} Let $(H', w)$ be the resulting bipartite signed edge-weighted  $t$-tree after adding $u$. Since $(H,w)$ and $(K, w)$ are $2k$-wide, we only need to show that any negative cycle of $\overline{\overline{(H',w)}}_{-2k}$ containing $u$ and a vertex which is not in $(K,w)$ has length at least $2k$. Assume to the contrary and let $C$ be a negative cycle in $\overline{\overline{(H',w)}}_{-2k}$ which contains $u$ and has length less than $2k$. By the structure of $(H',w)$, at least two vertices of $C$ (other than $u$) are in $K$. Let $v,v'$ be such the two vertices. Since $K$ is a $(t+1)$-clique, $v$ is adjacent to $v'$ in $(H',w)$. Thus there are two $(v-v')$-paths, say $R_1$ and $R_2$, in $\overline{\overline{(H',w)}}_{-2k}$ one of length $|w(vv')|$, the other of length $2k-|w(vv')|$. Moreover, one of $R_1$ and $R_2$ is negative and the other is positive. The cycle $C$ is also split into two $(v-v')$-paths, say $P$ and $Q$, whose total length is less than $2k$ and have opposite signs. Assuming $R_1$ and $P$ are of the opposite signs, we have two negative cycles, $R_1P$ and $R_2Q$, whose total length is less than $4k$. Thus at least one of the two has a length smaller than $2k$. But one of the cycles is in  $\overline{\overline{(H,w)}}_{-2k}$ and the other is in  $\overline{\overline{(K,w)}}_{-2k}$. However, both $(H, w)$ and $(K,w)$ are assumed to be $2k$-wide. This contradiction completes the proof.
	\end{proof}
	
	\subsection{Homomorphism of signed edge-weighted graphs}
	
	Given signed edge-weighted graphs $(G, w_1)$ and $(H, w_2)$, a homomorphism of $(G, w_1)$ to $(H, w_2)$ is to first apply a suitable switching on $(G, w_1)$ to get $(G, w'_1)$ and then map $G$ to $H$ in such a way that the weights of the edges are preserved with respect to $w'_1$ and $w_2$. 
	
	When certain conditions are met, then a mapping of $(G, w_1)$ to $(H, w_2)$ will be equivalent to a mapping of $\overline{\overline{(G,w_1)}}_{-2k}$ to $\overline{\overline{(H,w_2)}}_{-2k}$. The goal of the next subsection is to set such conditions.
	
	\subsection{Algebraic distances as weights}
	
	A key property that this work is built on is that when mapping a negative cycle to a negative cycle of the same length, the distances between pairs of vertices is preserved.  If the mapping is, furthermore, edge-sign preserving, then it preserves the algebraic distances. This is our motivation for defining signed edge-weighted graphs, so the algebraic distances will be the weights. Observe that the more constraint we have for a potential homomorphism, the more restricted the choices are. When the number of choices is limited to a very few, then it would be practical to either find a suitable mapping or verify that there is none. Moreover, in working with partial $t$-trees ($t\geq 2$) of a given negative girth, we may assume each edge is in a negative cycle of the said size. These observations lead to the following definitions and notations.

	\begin{definition} The \emph{signed $(G,\sigma)$-distance complete graph} is the signed edge-weighted complete graph $(K_{|V(G)|}, ad_{(G,\sigma)})$, where the weight of an edge $xy$ is the algebraic distance $ad_{(G,\sigma)}(x,y)$ in $(G,\sigma)$ between $x$ and $y$.
	\end{definition}
	
	Therefore, the edges of $(K_{|V(G)|}, ad_{(G,\sigma)})$ of weights 1 and $-1$ induce $(G,\sigma)$. 
	
	\begin{definition} Any spanning (weighted) subgraph of $(K_{|V(G)|}, ad_{(G,\sigma)})$ is called  a \emph{partial signed distance graph} of $(G,\sigma)$ or a \emph{partial signed $(G,\sigma)$-distance graph}. 
	\end{definition}

	Let $\widehat{B}$ be a signed bipartite graph of unbalanced girth $2k$. We define signed edge-weighted graph $(\widehat{B},w)_{-2k}$ as follows: for each pair  of vertices $x,y$ that lie on a common negative $2k$-cycle add an edge and assign the weight $w(xy)=ad_{\widehat{B}}(x,y)$. Thus, an original edge remains an edge (of weight $1$ or $-1$) if and only if it is in a negative $2k$-cycle. 
	
	\begin{lemma}\label{homomorphismOf(G,w)(B,w)} Let $(H,w_1)$ be a $2k$-wide signed edge-weighted graph satisfying $|w_1(e)|\leq k$ for every $e\in E(H)$ and let $\widehat{B}$ be a signed graph of unbalanced girth $2k$. Then any homomorphism, say $\phi$, of $\overline{\overline{(H,w_1)}}_{-2k}$ to $\widehat{B}$, induces a homomorphism of $(H,w_1)$ to $(\widehat{B},w)_{-2k}$. 
	\end{lemma}
	
	\begin{proof} 
		By renaming $w_1$ if necessary, we may assume $\phi$ is an edge-sign preserving homomorphism of $\overline{\overline{(H,w_1)}}_{-2k}$ to $\widehat{B}$.
		By the construction of $\overline{\overline{(H,w_1)}}_{-2k}$, if $xy$ is an edge of $H$ with $w(xy)=i$, then it is on a negative cycle of length $2k$ in $\overline{\overline{(H,w_1)}}_{-2k}$ in which the algebraic distance of $x$ and $y$ is $i$. Moreover, since $(H,w_1)$ is $2k$-wide, $i$ is the algebraic distance of $x$ and $y$ in $\overline{\overline{(H,w_1)}}_{-2k}$. Since $\widehat{B}$ is of unbalanced girth $2k$, the images $\phi(x)$ and $\phi(y)$, which must be on a negative $2k$-cycle, are also at algebraic distance $i$. This means that $\phi$ is an (edge-sign preserving) mapping of $(H,w_1)$ to $(\widehat{B},w)_{-2k}$. 
	\end{proof}
	
	\subsection{Unbalanced-girth and wideness}
	
	Let $(G,\sigma)$ be a signed bipartite partial $t$-tree of unbalanced-girth at least $2k$. Let $H$ be a $t$-tree completion of $G$ with a $(t+1)$-clique sequence $H_{t+1}, H_{t+2},\ldots,H_{t+l}$. That is, $V(G)=V(H)$ and $E(G)\subseteq E(H)$, $H$ is a $t$-tree. A weight function $w$ on the edges of $H$ is defined as follows:
	
	\begin{equation}\label{equ:weights}
		w(uv)=\left\{
		\begin{array}{ll}
			ad_{(G,\sigma)}(u,v),     &   {\rm if}\ d_{G}(u,v)\leq k;\\
			k,     &    {\rm if}\ d_{G}(u,v)>k{\rm\ and}\ d_{G}(u,v) {\rm\ has\ the\ same\ parity\ as}\ k;\\
			k-1,     &    {\rm if}\ d_{G}(u,v)>k{\rm\ and}\ d_{G}(u,v) {\rm\ has\ the\ same\ parity\ as}\ k-1.\\
		\end{array} \right.
	\end{equation}

	Since each weighted cycle of $(H, w)$ corresponds to a closed walk of the same parity in $(G, \sigma)$, the weighted length of every cycle in $(H, w)$ is even and thus $(H,w)$ is bipartite. Next, we show that every signed edge-weighted $(t+1)$-clique of $(H,w)$ is $2k$-wide.

	\begin{lemma}\label{lemma:2k-wide} Let $(G,\sigma)$ be  a signed bipartite partial $t$-tree of unbalanced-girth at least $2k$ and $H$ be a $t$-tree completion of $G$ with $w$ being the weight function on the edges of $H$ as defined in equation (1). Then every signed edge-weighted $(t+1)$-clique $(H_{t+i},w)$ of $(H, w)$ is $2k$-wide.
	\end{lemma}
	
	\begin{proof}
		By Theorem \ref{thm:complete weighted graph}, to prove that $(H_{t+i},w)$ is $2k$-wide for each $i\in \{1,2,\ldots,l\}$, it is enough to show that every signed edge-weighted triangle $(K_3,w)$ of $(H_{t+i},w)$ is $2k$-wide. Let $x,y,z$ be the three vertices of a triangle $T$ of $(H_{t+i},w)$.
		If $w(xy)w(xz)w(zy)<0$, then in $\overline{\overline{(T,w)}}_{-2k}$ a shortest negative cycle containing all three of $x$, $y$, and $z$ is of length $|w(xy)|+|w(xz)|+|w(zy)|$. By the definition of $w$, either this corresponds to the length of a negative closed walk of $(G, \sigma)$, which then is at least $2k$, or at least one of $|w(xy)|$, $|w(xz)|$,  $|w(zy)|$ is larger than $k$, but then the sum of the other two is also at least $k$ by the triangular inequality.

		If $w(xy)w(xz)w(zy)>0$, then assuming $|w(xy)|\geq |w(xz)| \geq |w(zy)|$, in $\overline{\overline{(T,w)}}_{-2k}$ a shortest negative cycle containing all three of $x$, $y$, and $z$ is of length $2k-|w(xy)|+|w(xz)|+|w(zy)|$. But then the triangular inequality implies that this length is at least $2k$.
	\end{proof}



	
	\subsection{List of $2k$-wide $(t+1)$-cliques and a $(t,2k)$-closed set}

	The two key notions for the main characterization of bounds are the following definition, noting that it is essential not to consider a labeling of vertices for the elements of the list (in \Cref{def:ListOf2kWide}), while labeling of vertices is essential for a set of cliques to be $(t,2k)$-closed (in \Cref{t2k-closed}).
	
	\begin{definition}\label{def:ListOf2kWide} List of $2k$-wide $(t+1)$-cliques, denoted $\mathcal{L}(t+1,2k)$, is the set of all $2k$-wide $(t+1)$-cliques without labeling on the vertices of the cliques.
	\end{definition}
	
	Noting that the label of each edge is a (nonzero) number between $-k$ and $k$, we have a trivial upper bound of $(2k)^{t+1 \choose 2}$ on the number of elements of $\mathcal{L}(t+1,2k)$. 
	
	\begin{definition}\label{t2k-closed} Given a signed edge-weighted $2k$-wide graph $(H, w)$, a set $\mathcal{W}$ of $(t+1)$-cliques is said to be $(t,2k)$-closed if the following holds: for any signed edge-weighted $(t+1)$-clique $(K,w)\in \mathcal{W}$, with vertices labeled $v_1,v_2,\ldots,v_t,v_{t+1}$, if $X\in \mathcal{L}(t+1,2k)$ is such that the subgraph induced by $K-v_i$ is isomorphic to $X-x$ for a vertex $x$ of $X$, with $\xi$ being the isomorphism, then there is a clique $(K', w)\in \mathcal{W}$ which, first of all, contains all vertices in $K-v_i$, secondly, $\xi$ extends to an isomorphism of $(K',w)$ and $X$. 
	\end{definition}
	
	The condition that there should be a $(K',w)$ for every isomorphism $\xi$ implies that any $(t, 2k)$-closed set of $(t+1)$-cliques must have significantly many more elements than the size of the list $\mathcal{L}(t+1,2k)$. In the next lemma we prove the bare minimum though; that in any $(t, 2k)$-closed set there is at least one clique isomorphic to any given element of $\mathcal{L}(t+1,2k)$.

	\begin{lemma}\label{isomorphic copy} Let $t$ and $k$ be positive integers, $(H, w)$ be a $2k$-wide signed edge-weighted graph, and $\mathcal{W}$ be a nonempty set of $(t+1)$-cliques of $(H, w)$ which is $(t,2k)$-closed. 
	Then for any $(t+1)$-clique $X\in \mathcal{L}(t+1,2k)$, there is an isomorphic copy $W$ of $X$ in $\mathcal{W}$.
	\end{lemma}

	\begin{proof} Let $a$ be the even value among $k$ and $k-1$ and let $X_k$ be the element of $\mathcal{L}(t+1,2k)$ where the weight of each edge is  $a$.
		Since $\mathcal{W}$ is not empty, there is a $(t+1)$-clique $W$ in $\mathcal{W}$. Let $v$ be a vertex of $W$ which is incident to at least one edge of weight different from $a$, let $W_v$ be  one of the two elements $X^0$ or $X^1$ of $\mathcal{L}(t+1,2k)$ defined as follows. 
		
		In $X^0$ all weights are as in $W$ except for weights of the edges incident to $v$. For each edge incident with $v$ if its weight (in $W$) is even, then we replace it with $k$ and if the weight is odd, then replace it with $k-1$. To define $X^1$ we do the same except we change to role of odd and even. 
		  
		It can be checked readily that both $X^0$ and $X^1$ are members of $\mathcal{L}(t+1,2k)$. Since $\mathcal{C}$ is $(t, 2k)$-closed, there are $(t+1)$-cliques in $\mathcal{W}$ one isomorphic to $X^0$ and the other isomorphic to $X^1$. Of the two, the one with more edges of weight $a$ is the one we take to be $W_v$. 
		
		By repeating this process, starting from any element $W$ of $\mathcal{W}$ and using the assumption that $\mathcal{W}$ is $(t, 2k)$-closed, we get to an element of $\mathcal{W}$ which is isomorphic to $X_k$. As the process of the changing vertices is reversible, continuing from this particular element of $\mathcal{W}$ we can get an isomorphic copy of any element in $\mathcal{L}(t+1,2k)$.		
	\end{proof}

	\section{Bounding signed bipartite partial $t$-trees}
	
	We now have all the tools to state a necessary and sufficient condition for a signed bipartite graph of negative girth $2k$ to bound the class  $\mathcal{SBPT}_{t,2k}$ of signed bipartite  partial $t$-trees of unbalanced girth at least $2k$. As it will be pointed out, the theorem can be used to practically determine if a given signed bipartite graph bounds $\mathcal{SBPT}_{t,2k}$.

	\begin{theorem} \label{thm:NecessarySufficient}  A signed bipartite graph $\widehat{B}$ of unbalanced girth $2k$ bounds $\mathcal{SBPT}_{t,2k}$ if and only if there exists a partial $\widehat{B}$-distance graph $(\widehat{B},w)$ with a non-empty set $\mathcal{W}$ of $(t+1)$-cliques which is $(t,2k)$-closed.
	\end{theorem}

	\begin{proof} To prove that the condition is sufficient, we consider a signed graph $(G,\sigma)\in \mathcal{SBPT}_{t,2k}$. Let $H$ be a $t$-tree completion of $G$ with vertex sequence $v_1,v_2,\ldots,v_n$ and $(t+1)$-clique sequence $H_{t+1}, H_{t+2},\ldots,H_{n}$. Let $(H,w')$ be the signed edge-weighted graph, where $w'$ is the weight function defined in \Cref{equ:weights}. Our goal is to find a homomorphism of $(H,w')$ to $(\widehat{B},w)$ noting that this will also be a homomorphism of $(G,\sigma)$ to $\widehat{B}$. 
		
		By Lemma \ref{lemma:2k-wide}, we know that $(H_{t+1},w')$ is $2k$-wide. By Lemma \ref{isomorphic copy}, we can find an isomorphic copy of $(H_{t+1},w')$ in $\mathcal{W}$. Let $u_1,u_2,\ldots,u_{t+1}$  be the images of $v_1,v_2,\ldots,v_{t+1}$ in $(\widehat{B},w)$, respectively. Consider the second clique $(H_{t+2},w')$. We know that $|V(H_{t+1})\cap V(H_{t+2})|=t$. Without loss of generality, assume that $V(H_{t+1})\cap V(H_{t+2})=\{v_1,v_2,\ldots,v_t\}$. For $i\in \{1,2,\ldots,t\}$, let $a_i=w'(v_iv_{t+2})$. Since $\mathcal{W}$ is $(t,2k)$-closed, there exists $u_{t+2}\in V(B)$ such that $w(u_iu_{t+2})=a_i$ and that $u_1,u_2, \ldots, u_t, u_{t+2}$ induces a clique in $\mathcal{W}$. We map $v_{t+2}$ to $u_{t+2}$. Continuing this process, we map $(H,w')$ to $(B,w)$.
		
		For the converse statement, let $\widehat{B}$ be a signed bipartite graph of unbalanced girth $2k$ which bounds $\mathcal{SBPT}_{t,2k}$. Let $(\widehat{B},w)$ be a $\widehat{B}$-partial distance graph whose edges are all pairs of vertices that belong to a common $C_{-2k}$ where for any such pair $x$ and $y$, the weight $w(xy)$ is the algebraic-distance of $x$ and $y$ in the common $C_{-2k}$. Observe that, since $\widehat{B}$ is of unbalanced girth $2k$, the value $w(xy)$ is well-defined and is the same as the algebraic distance of $x$ and $y$ in $\widehat{B}$. 
		
		For a given signed bipartite  partial $t$-tree $(G,\sigma)$, let $(H,w')$ be as defined by \Cref{equ:weights}. Observe that, for $t\geq 2$, $\overline{\overline{(H,w')}}_{-2k}$ is also a signed bipartite partial $t$-tree of unbalanced girth $2k$. Hence, by our assumption, it also maps to $\widehat{B}$. Any such a mapping then induces a mapping of $(H, w')$ to $(\widehat{B},w)$. Thus we may assume that all $2k$-wide bipartite signed edge-weighted $t$-trees map to $(\widehat{B},w)$.      
			
	    Let $\mathcal{W}$ be a minimal set of $(t+1)$-cliques in $(\widehat{B},w)$ such that for every bipartite signed $2k$-wide $t$-tree $(H,w')$, there is a homomorphism of $(H,w')$ to $(\widehat{B},w)$ satisfying that the image of every $(t+1)$-clique of $(H,w')$ is in $\mathcal{W}$. 
		
		The assumption on minimality implies the following: for every (weighted) $(t+1)$-clique, say $K_{t+1}^i$ of $\mathcal{W}$, there exists a bipartite signed $2k$-wide $t$-tree $(H_i,w_i')$ such that any homomorphism of $(H_i,w_i')$ to $(\widehat{B},w)$ which maps all $(t+1)$-cliques of $(H_i,w_i')$ to $\mathcal{W}$ must map at least one of such cliques to $K_{t+1}^i$.
		
		We use this property to show that the set $\mathcal{W}$ is $(t,2k)$-closed. To this end first observe that $\mathcal{W}$ is not empty because any weighted $(t+1)$-clique which is $2k$-wide must have an image in $\mathcal{W}$. Let $W_1$ be a clique in $\mathcal{W}$ with vertices $v_{i_1}, v_{i_2}, \ldots, v_{i_{t+1}}$ and assume removing $v_{i_{t+1}}$ and adding a vertex $x$ which is adjacent to $v_{i_{j}}$ by an edge of weight $w(xv_{i_j})$ we get a $2k$-wide clique. That is to say, the unlabeled clique is in the list  $\mathcal{L}(t+1,2k)$.
		
		Our goal is to show that there is a clique $W'$ in $\mathcal{W}$ on vertices $v_{i_1}, v_{i_2}, \ldots, v_{i_{t}}, v_x$ such that each edge $v_{i_{j}}v_x$, $j=1, \ldots, t$, has the same weight as $w(xv_{i_j})$.
		
		Toward this consider $(H_i,w_i')$ and for each isomorphic copy $P$ of $W_1$, add a vertex $x_P$ and join it to the $t$ vertices of $P$ that correspond to $v_{i_1}, v_{i_2}, \ldots, v_{i_{t}}$ under the assumed isomorphism. Thus for each such a clique $W_1$ the number of $x_p$ that are added is the number of automorphisms of $W_1$.   
		
		Let $H_i^*$ be the resulting signed weighted $t$-tree. By  Lemma~\ref{thm:a new vertex}, $H_i^*$ is also $2k$-wide. Thus it admits a homomorphism $\varphi$ to $(\widehat{B},w)$ such that each $(t+1)$-clique is mapped to an element of $\mathcal{W}$. Observe that $\varphi$ also maps $(H_i,w_i')$ to $(\widehat{B},w)$ such that all $(t+1)$-cliques are mapped to the elements of $\mathcal{W}$. Thus, by our assumption about $(H_i,w_i')$, at least one of the $(t+1)$-cliques, say $P$, is mapped to $W_1$. Considering the isomorphism of $P$ to $W_1$ induced by $\varphi$, we have an associated vertex $x_P$, and $\varphi(x_P)$ is the vertex $v_x$ that we are looking for.		
	\end{proof}
	
	\section{Bounding signed partial 3-trees}\label{sec:PCbounds3-trees}
	
	As one of the applications of \Cref{thm:NecessarySufficient}, we show that signed projective cube of dimension $2k-1$, which is a signed bipartite graph of negative girth $2k$, bounds the class of signed bipartite partial 3-trees of negative girth at least $2k$. 
	This is related to several conjectures in (potential) extensions of the 4-color theorem and provides support for such conjectures. We first recall the definition of signed projective cube. For a comprehensive study of these signed graphs we refer to \cite{CNS25}. 
	
	\begin{definition} The \emph{signed projective cube} of dimension $n$, denoted $SPC(n)$, is the signed graph on vertex set $\mathbb{Z}_{2}^{n}$ where two vertices at the hamming distance 1 are adjacent with a positive edge and those at hamming distance $n$ are adjacent with a negative edge. 
	\end{definition}
	
	Thus we may label each positive edge with an element $e_i$ of the standard basis of $\mathbb{Z}_{2}^{n}$ and each negative edge with the element $J=11\cdots1$. 
	Let $S_n=\{e_1, \ldots, e_n, J\}$. Observe that the only linear relations among the elements of $S_n$ are $2s=0$ for any element, and $e_{1}+e_{2}+\cdots+e_{n}+J=0$. As sum of the weights of each cycle is $0$, there are two types of cycles in $SPC(n)$: 
	\begin{itemize}
		\item Positive cycles each of which contains each element of $S_n$ an even number of times (0 being allowed). 
		\item Negative cycles each of which contains each element of $S_n$ an odd number of times, thus in particular having length at least $n+1$.
	\end{itemize}
	
	An equivalent definition of $SPC(n)$ is as follows. Vertices of $SPC(n)$ are the (unordered) pairs of complementary subsets of $S_n$. Two vertices $\{A, \bar{A}\}$ and $\{B, \bar{B}\}$ are adjacent if $A$ contains $B$ and has exactly one more element. The sign of such edge is negative only if the extra element of $A$ is $J$. 
	
	In view of this definition, the distance of two vertices $u=\{A, \bar{A}\}$ and $v=\{B, \bar{B}\}$ is the order of the smaller of the two symmetric differences: $A\oplus B$ and $A \oplus \bar{B}$. Observing that these two sets are complementary subsets of $S_{n}$, the one which does not contain $J$ is denoted by $S^+(u,v)$ and the one that contains $J$ is denoted by $S^-(u,v)$. The sign of the algebraic distance of $u$ and $v$ then is based on whether $|S^+(u,v)|\leq k$ or $|S^-(u,v)|<k$, the former being associated with the positive sign.  The complementary sets $S^+(u,v)$ and $S^-(u,v)$ associated to each pair of vertices present the relative position of vertices $u$ and $v$ on any $C_{-(n+1)}$ that they both belong to; the edges of the negative $u-v$ path of any such cycle are labeled with (distinct) elements of $S^-(u,v)$ and the edges of the positive $u-v$ path are labeled with (distinct) elements of $S^+(u,v)$.

	Thus for each set of three (distinct) vertices, $x$, $y$, and $z$, there are essentially four minimal negative closed walks passing through the three of them and four positive ones. The negative ones correspond to: $\{S^-(x,y), S^-(y,z), S^-(z,x)\}$, $\{S^-(x,y), S^+(y,z), S^+(z,x)\}$, $\{S^+(x,y), S^-(y,z), S^+(z,x)\}$, and $\{S^+(x,y), S^+(y,z), S^-(z,x)\}$, while the positive ones correspond to the other four. 
	
	Recall that for each of these triplets representing minimal negative walks through $x$, $y$ and $z$, each element of $S_{n}$ appears in an odd number of them. That is to say, either exactly in one of them or in all three. The number of elements that appear in all three of them, say in $S^-(x,y)\cap S^-(y,z)\cap S^-(z,x)$ then is only a function of the size of these sets. To be precise, that is, $\tau_{xyz}=\frac{|S^-(x,y)|+|S^-(y,z)|+|S^-(z,x)|-|S_n|}{2}$. It follows from this formula that the values of $|S^-(x,y)|, |S^-(y,z)|$, and $|S^-(z,x)|$ determines the types of triplet of vertices. This is precised in the following theorem.

		\begin{theorem}\label{thm:PC(2k)TripleTransitive}
		For any positive integer $k$, the signed projective cube $SPC(k)$ is triple-transitive. That is to say, given two triplet of vertices, $\{x,y,z\}$ and $\{u,v,t\}$, if after a switching we have that $(|S^-(x,y)|, |S^-(y,z)|,$ $|S^-(z,x)|)=(|S^-(u,v)|, |S^-(v,t)|, |S^-(t,u)|)$, then an automorphism of $SPC(k)$ maps $(x,y,z)$ to $(u,v,t)$.
	\end{theorem}

	\begin{proof}
		The assumption implies that, subject to some switching, $\tau_{xyz}=\tau_{uvt}$. Thus each triplet partitions $S_{n}$ into four parts:  $S_n=S_{xyz}\cupdot S'_{xy}\cupdot S'_{yz} \cupdot S'_{zx}$ and $S_n=S_{uvt}\cupdot S'_{uv}\cupdot S'_{vt} \cupdot S'_{tu}$ where $S_{xyz}=S^-(x,y)\cap S^-(y,z)\cap S^-(z,x)$, $S'_{xy}=S_{xy}-S_{xyz}$, $S'_{yz}=S_{yz}-S_{xyz}$, and $S'_{zx}=S_{zx}-S_{xyz}$, with $S_{uvt}, S'_{uv}, S'_{vt},$ and $S'_{tu}$ being defined similarly.		
		By our assumption we have $|S_{xyz}|=|S_{uvt}|=\tau_{xyz}=\tau_{uvt}$, $|S'_{xy}|=|S'_{uv}|=|S_{xy}|-\tau_{xyz}=|S_{uv}|-\tau_{uvt}$, $|S'_{yz}|=|S'_{vt}|=|S_{yz}|-\tau_{xyz}=|S_{vt}|-\tau_{uvt}$, and $|S'_{zx}|=|S'_{tu}|=|S_{zx}|-\tau_{xyz}=|S_{tu}|-\tau_{uvt}$. 
		
		Let $\mu$ be a permutation of $S_n$ which maps $S_{xyz}$, $S'_{xy}$, $S'_{yz}$, and $S'_{zx}$ respectively to $S_{uvt}$, $S'_{uv}$, $S'_{vt}$, and $S'_{tu}$. 		
		An automorphism of $SPC(k)$ which maps $(x,y,z)$ to $(u,v,t)$ is a composition of the following three basic automorphisms: 1. a switching that is preformed in order to have $(|S^-(x,y)|, |S^-(y,z)|, |S^-(z,x)|)=(|S^-(u,v)|, |S^-(v,t)|, |S^-(t,u)|)$, 2. the permutation $\mu$ of $S_n$, 3. the mapping $\phi(a)=a+(x+u)$ where additions are taken in $\mathbb{Z}_2^k$. The last item maps $x$ to $u$, the first two then ensure that $y$ is mapped to $v$ and $z$ is mapped to $t$. We leave the verification of fine details to the reader.
	\end{proof}
	
	Recall that signed $SPC(n)$-distance complete graph,  $(K_{SPC(n)}, ad)$, is the signed edge-weighted complete graph whose vertices are the elements of $\mathbb{Z}_2^n$ where the weight of an edge $uv$ is the algebraic distance of $u$ and $v$ in $SPC(n)$, noting that the distance of two vertices is at most $\lfloor \frac{n}{2} \rfloor$. Each edge $uv$ of this signed edge-weighted graph is furthermore labeled with the partition $S_{uv}\cupdot \bar{S}_{uv}$ of $S_n$, noting that the size of the smaller part is the distance of $uv$ and the sign of the corresponding weight depends on whether the smaller part contains $J$ or not. It is proved in \cite{CN20} that $PC(2k)$ bounds the class of partial 3-trees of odd-girth at least $2k+1$. The following theorem implies analogue of this result for signed bipartite partial 3-trees.

	\begin{theorem}\label{(3,2k)closed}
		The set $\cal{W}$ consisting of all weighted 4-cliques of $(K_{SPC(2k-1)},ad)$ is $(3,2k)$-closed.
	\end{theorem}
	
	\begin{proof} 
		Since $(K_{SPC(2k-1)},ad)$ is triple-transitive, it would be enough to show that for each element $L$ of the list $\mathcal{L}(4, 2k)$ there is a $(K_4,w)$ in $(K_{SPC(2k-1)},ad)$ realizing it. That is because given a triple $x,y,z$ of the vertices if we need to find a vertex at distance $\alpha,\beta, \gamma$ from $x,y,z$, we consider $(K_4,w)$ whose weights are the same as the weights of the triangle and $a,b,c$. We then map the triangle to $x,y,z$ by an automorphism. The image of the fourth vertex is what are looking for.  
		
		So we consider a $2k$-wide $(K_4,w)$, as labeled in \Cref{fig:K4Exist}. For the rest of this proof, noting that each edge is labeled with a pair $(a, a')$ of integers such that $aa'<0$, and $|a|+|a'|=2k$, to represent the edge we select the negative value among $a$ and $a'$. What remains to prove is to assign a $|w(e)|$-subset of $S_{2k-1}$ to each edge $e$ of the given $(K_4, w)$ such that each of them contains $J$, and the sum (in $\mathbb{Z}_2^{2k-1}$) of the weights of all the edges of any of the cycles is $0$. Having done so, to realize the given $(K_4,w)$ in  $(K_{SPC(2k-1)},ad)$ we take any vertex of it as $x$, then $y$, $z$,  $v$ are respectively $x+\sum\limits_{r\in S_{xy}} r$, $x+\sum\limits_{r\in S_{xz}} r$, and $x+\sum\limits_{r\in S_{xv}} r$.
		
		Let $t_0=\frac{a+b+c-2k}{2}$, $t_1=\frac{a+\beta+\gamma-2k}{2}$, $t_2=\frac{\alpha + b+\gamma-2k}{2}$, and $t_3=\frac{\alpha +\beta+c-2k}{2}$. We know that $t_0, t_1, t_2, t_3 \geq 1$. Without loss of generality, assume that $t_0\leq t_1\leq t_2\leq t_3$. We claim that $(t_1-t_0)+(t_2-t_0)+(t_3-t_0)\leq 2k-t_0$. That is because in   $\overline{\overline{(K_4,w)}}_{-2k}$ we have a negative cycle of length $a+(2k-\beta)+ (2k-\gamma)$. Since $(K_4,w)$ is $2k$-wide, we have $a+2k-\beta+2k-\gamma\geq 2k$, equivalently, $a\geq \frac{a+\beta+\gamma-2k}{2}=t_1$. Similarly, we have $b \geq t_2$ and $c \geq t_3$. We get $t_1+t_2+t_3\leq a+b+c=2k+2t_0$.

        We now partition $S_{2k-1}$ into seven subsets: $S_{t_0}, S_{t_1-t_0}, S_{t_2-t_0}, S_{t_3-t_0},$ $S_{a-t_1}, S_{b-t_2}, S_{c-t_3}$.
         Some of these subsets could potentially be empty. The index of each subset presents the size of the set.
         
        To complete the assignment then we choose $J$ and $t_0-1$ elements of $S_{2k-1}$ and assign it to all $S_{e}$. The following then is an assignment of subsets of $S_{2k-1}$ to the edges which satisfies the required conditions: 
        
        \begin{itemize}
        	\item $S_a=S_{t_0}\cup S_{t_1-t_0}\cup S_{a-t_1}$,
        	\item $S_b=S_{t_0}\cup S_{t_2-t_0}\cup S_{b-t_2}$,
        	\item $S_c=S_{t_0}\cup S_{t_3-t_0}\cup S_{c-t_3}$,
        	\item $S_{\alpha}=S_{t_0}\cup S_{t_2-t_0}\cup S_{t_3-t_0} \cup S_{a-t_1}$,
        	\item $S_{\beta}=S_{t_0}\cup S_{t_1-t_0}\cup S_{t_3-t_0} \cup S_{b-t_2}$,
        	\item $S_{\gamma}=S_{t_0}\cup S_{t_1-t_0}\cup S_{t_2-t_0} \cup S_{c-t_3}$.
        \end{itemize}

	\begin{figure}[ht]
		\centering 
		\begin{minipage}{0.45\textwidth}
			\centering
			\begin{tikzpicture}
				\foreach \i in {1,2,3}
				{
					\draw[rotate=-120*(\i)] (0, 3) node[circle, fill=white, draw=black!80, inner sep=.7mm, minimum size=1mm]  (\i){};
					\draw (0, 0.1) node[circle, fill=white, draw=black!80, inner sep=.7mm, minimum size=1mm]  (4){};
					\node [left, rotate=60] at (-1.3, 1) {$a$};
					\node [right,  rotate=-60] at (1.4, 1) {$b$};
					\node [right,  rotate=-30] at (.7, -0.15) {$\alpha$};
					\node [left,  rotate=30] at (-0.7, -0.15) {$\beta$};
					\node [below] at (0, -1.5) {$c$};
					\node [right] at (0, 1.5) {$\gamma$};
					
					\node [below] at (0, 0) {$v$};
					\node [right] at (2.7, -1.5) {$x$};
					\node [left] at (-2.7, -1.5) {$y$};
					\node [above] at (0, 3.1) {$z$};
				}
				
				\foreach \i/\j in {1/2, 1/3, 2/3, 1/4, 2/4, 3/4}
				{
					
					\draw[thick] (\i) -- (\j);
				}
				
			\end{tikzpicture}
			
		\end{minipage}
		\begin{minipage}{0.45\textwidth}
			\centering
			\begin{tikzpicture}
				\foreach \i in {1,2,3}
				{
					\draw[rotate=-120*(\i)] (0, 3) node[circle, fill=white, draw=black!80, inner sep=.7mm, minimum size=1mm]  (\i){};
					\draw (0, 0.1) node[circle, fill=white, draw=black!80, inner sep=.7mm, minimum size=1mm]  (4){};
					\node [left, rotate=60] at (-1, 2) {$S_{t_0}\cup S_{t_1-t_0}\cup S_{a-t_1}$};
					\node [right, rotate=-60] at (1.2, 1.8) {$S_{t_0}\cup S_{t_2-t_0}\cup S_{b-t_2}$};
					\node [right,  rotate=-30] at (.8, -0.12) {$S_{\alpha}$};
					\node [left,  rotate=30] at (-0.8, -0.12) {$S_{\beta}$};
					\node [below] at (0, -1.5) {$S_{t_0}\cup S_{t_3-t_0}\cup S_{c-t_3}$};
					\node [right] at (0, 1.5) {$S_{\gamma}$};
				}
				
				\foreach \i/\j in {1/2, 1/3, 2/3, 1/4, 2/4, 3/4}
				{
					
					\draw[thick] (\i) -- (\j);
				}

			\end{tikzpicture}

		\end{minipage}
		\caption{Finding $(K_4,w)$ in $(K_{SPC(2k-1)},ad)$.}
		\label{fig:K4Exist}
		
	\end{figure}
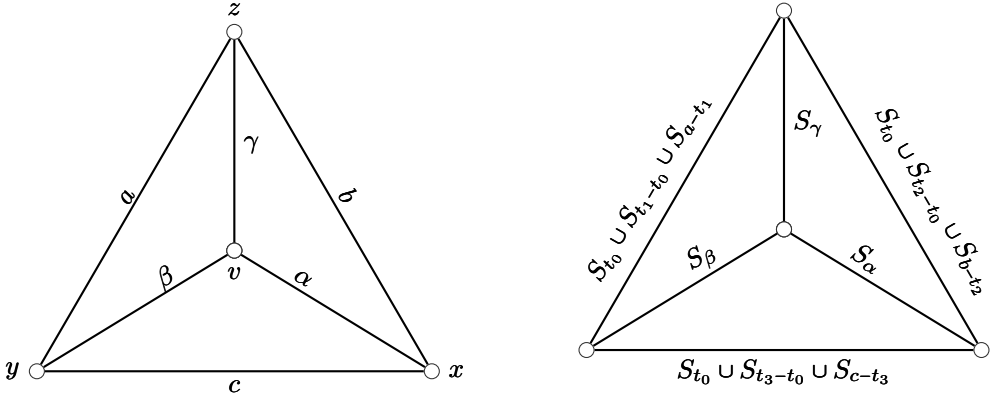
		\end{proof}
		
	
	We remark that in this proof, we chose the elements of $S_{t_0}$ to be contained in all the six sets. For some choices of weights, this is the only possibility. However, when $k$ is big enough, for most choices of the weights there are other possibilities of assigning these sets. We noted earlier that signed projective cubes are triple transitive. This observation implies that this is the limit of the property and that there are sets of four vertices that induces the same distance system, but they cannot be mapped to each other by an automorphism.  
	
	Applying  \Cref{thm:NecessarySufficient} we conclude that $SPC(2k-1)$ bounds the class of signed bipartite partial 3-trees of unbalanced girth at least $2k$. We may combine this result with the result of \cite{CN20} for a uniform presentation based on the following notation. 
	
	Given a positive integer $k$, let $\mathcal{O}_k$ be the set of signed cycles that do not map to $C_{-k}$ (the negative cycle of length $k$). Thus $\mathcal{O}_k$ consists of positive odd cycles, negative cycles of length less than $k$, and negative cycles whose parity of length is different from that of $k$. The class $\mathcal{SPT}_{t, -k}$ then is the class of signed partial $t$-trees which do not contain any element of $\mathcal{O}_k$ as an (induced) subgraph. Observe that $SPC(k-1)$ also does not have an element of $\mathcal{O}_k$ as a subgraph.
	
	\begin{theorem}\label{thm:Partial3Tree+Ok}
		A signed partial 3-tree $(G, \sigma)$ maps to $SPC(k)$ if and only if it has no element of $\mathcal{O}_{k+1}$ as a subgraph.
	\end{theorem} 
	
	The result supports a conjecture of the second author and B. Guenin in extension of the 4-color theorem. The conjecture claims that the class of partial 3-trees in this theorem can be replaced with the class of $K_5$-minor free graphs. A deeper extension based on the notion of minor for signed graphs is the following:

	\begin{conjecture}\label{conj:NoOddK5ToSPC}
		Any signed graph $(G, \sigma)$ which has neither a $(K_5,-)$-minor nor induces an element of $\mathcal{O}_{k+1}$ maps to $SPC(k)$.
	\end{conjecture} 

    We recall that a signed graph $(H, \pi)$ is a minor of $(G, \sigma)$ if it can be obtained from $(G,\sigma)$ by a sequence of: 1. deletion (vertices or edges), 2. contraction of positive edge, 3. switching.
    
   \subsection{Connection to Extended Double Cover operation}
	A deeper look at this conjecture suggests what makes $SPC(k)$ special is their inductive definition based on the notion of Extended Double Cover: Given a signed graph $\widehat{B}$, its extended double cover, denoted $EDC(\widehat{B})$, has as its vertices two copies  $V^+(\widehat{B})$, and $V^{-}(\widehat{B})$ of the vertices of $\widehat{B}$ where $x^{+}x^{-}$ are the (only) negative edges of $\widehat{B}$, if $xy$ is a positive edge of $\widehat{B}$, then both  $x^{+}y^{+}$ and $x^{-}y^{-}$  are positive edges, and if $xy$ is a negative edge of $\widehat{B}$, then both $x^{+}y^{-}$ and $x^{-}y^{+}$ are positive edges.
	
	The following is essentially Theorem 28 of \cite{CNS25}, see also Section 7.2 of \cite{CNS25}.

	\begin{theorem}
	Let $\mathcal{C}$ be a minor closed family of graphs satisfying the following: for every $i\leq k$, the class of signed graphs $\mathcal{SC}_i=\{(G, \sigma) \mid G \in \mathcal{C}, (G, \sigma) \text{ contains no element of } \mathcal{O}_{i}\}$ is bounded by $SPC(i-1)$. If the class $\mathcal{SC}_l$, $l<k$, is bounded by a signed graph $\widehat{B}$, then the class $\mathcal{SC}_{l+1}$ is bounded by $EDC(\widehat{B})$. 
    \end{theorem}
    
    Taking $\mathcal{C}$ to be the class of partial $3$-trees, since the key assumption of the theorem holds for all values of $k$, each time we have a homomorphism bound $\widehat{B}$ for the class $\mathcal{SPT}_{3,-k}$, the signed graph $EDC(\widehat{B})$ is a bound for  $\mathcal{SPT}_{3, -(k+1)}$. Some special mentions are as follows.
    
    Let $\mathcal{C}$ be the class of partial 3-trees that are also planar. Thus $\mathcal{C}$ consists of graphs with five forbidden minors: the four graphs of \Cref{fig:forbiddenminorsofpartial3trees} to ensure being a partial 3-tree, with an addition of $K_{3,3}$ to ensure being planar. 
    
    Restricting a result of \cite{NPW22} from signed planar graphs to signed planar  graphs that are also partial 3-trees, we have: 
    
    \begin{theorem}
    	The class $\mathcal{SC}_8$ consisting of signed bipartite planar graphs which are also partial 3-trees and induce no element of $\mathcal{O}_8$ as subgraph are bounded by $C_{-4}$. 
    \end{theorem}

    Noting that $EDC(C_{-4})$ is switching equivalent to the signed graph on circular clique $C(8,3)$ where all edges are negative and that each signed graph with no induced subgraph from $\mathcal{O}_8$ can be switched to have no positive edge, we have the following.


    \begin{corollary}
    	Every planar graph which is also a partial $3$-tree and has odd girth at least 9 admits a circular $\frac{8}{3}$-coloring. 
    \end{corollary}
	
	A special case of Jaeger-Zhang conjecture states that every planar graph of odd girth 9 maps to $C_5$ which is the circular clique $C(5,2)$. If the conjecture holds, it would then not be surprise that it would hold for the larger class of $K_5$-minor-free graphs of odd-girth $9$, and thus for partial $3$-trees of odd girth at least 9. Hence, it is of special interest to replace $\frac{8}{3}$ of this corollary with the best possible value. It should be noted that for the smaller class of partial 2-trees, the best upper bound for the circular chromatic number based on the odd girth is given in \cite{PZ02}.

	\section{Application to edge-coloring}\label{sec:EdgeColoring}
	
    In this section, we apply our result to support Seymour's conjecture on edge-coloring planar multigraphs.
	
	If an $r$-regular multigraph $G$ is $r$-edge colorable, then each color class is a perfect matching. If $X$ is an odd set of vertices, then, being a perfect matching, each color class has at least one edge in the edge-cut $(X,X^c)$ which is referred to as an odd cut. Hence, in an $r$-edge-colorable $r$-regular graph there are at least $r$ edges in any odd edge-cut. This necessary condition is generally far from being sufficient. Petersen graph has its fame for being an example (or rather counterexample) for the case $r=3$ of this.
	
	In 1880, while Kempe's proof for the four-color theorem was believed, Tait introduced a reformulation of the four-color conjecture claiming that every bridge-less cubic planar graph is 3-edge-colorable.

	Extending this reformulation of the four-color theorem, about a century later, P.D. Seymour \cite{S79} conjectured that with an added assumption of planarity, the above mentioned necessary condition for $r$-edge-coloring of $r$-regular multigraph is also a sufficient one.
	
	\begin{conjecture} Any $r$-regular planar multigraph with no odd cut of size smaller than $r$ is $r$-edge colorable.
	\end{conjecture}
	
\begin{figure}[ht]
		\centering 
		
		\begin{minipage}{0.24\textwidth}
			\centering
			\begin{tikzpicture}
				\foreach \i in {1,...,5}
				{
					\draw[rotate=-72*(\i)] (0, 1.5) node[circle, fill=white, draw=black!80, inner sep=.7mm, minimum size=1mm]  (\i){};
				}
				
				\foreach \i/\j in {1/2, 1/3, 1/4, 1/5, 2/3, 2/4, 2/5, 3/4, 3/5, 4/5, 5/1}
				{
					
					\draw[thick] (\i) -- (\j);
				}

			\end{tikzpicture}
			
		\end{minipage}
		\begin{minipage}{0.24\textwidth}
			\centering
			\begin{tikzpicture}
				
				\foreach \i in {1,2,3}
				{
					\draw (0, 0.7) node[circle, fill=white, draw=black!80, inner sep=0.7mm, minimum size=1mm]  (1){};
					\draw (0, 2) node[circle, fill=white, draw=black!80, inner sep=0.7mm, minimum size=1mm]  (2){};
					\draw (0, 3.3) node[circle, fill=white, draw=black!80, inner sep=0.7mm, minimum size=1mm]  (3){};
				}
				
				\foreach \i in {4,5,6}
				{
					\draw (2, 0.7) node[circle, fill=white, draw=black!80, inner sep=0.7mm, minimum size=1mm]  (4){};
					\draw (2, 2) node[circle, fill=white, draw=black!80, inner sep=0.7mm, minimum size=1mm]  (5){};
					\draw (2, 3.3) node[circle, fill=white, draw=black!80, inner sep=0.7mm, minimum size=1mm]  (6){};
				}
				
				\foreach \i/\j in {1/4, 2/5, 3/6, 1/5, 1/6, 2/4, 2/6, 3/4, 3/5}
				{
					
					\draw[thick] (\i) -- (\j);
				}

			\end{tikzpicture}
			
		\end{minipage}
		\begin{minipage}{0.24\textwidth}
			
			\centering
			\begin{tikzpicture}
				\foreach \i in {1,...,5}
				{
					\draw[rotate=-72*(\i)] (0, 1.1) node[circle, fill=white, draw=black!80, inner sep=0.7mm, minimum size=1mm]  (\i){};
					
					\draw (0, 2) node[circle, fill=white, draw=black!80, inner sep=0.7mm, minimum size=1mm]  (6){};
					
					\draw (0, 0) node[circle, fill=white, draw=black!80, inner sep=0.7mm, minimum size=1mm]  (7){};
				}
				
				\foreach \i/\j in {1/2, 2/3, 3/4, 4/5, 5/1, 1/7, 2/7, 3/7, 4/7, 5/7, 5/6}
				{
					
					\draw[thick] (\i) -- (\j);
					\draw[bend left=20, thick] (6) to (1);
					\draw[bend right=20, thick] (6) to (4);
					\draw[bend left=80, thick] (6) to (2);
					\draw[bend right=80, thick] (6) to (3);
				}

			\end{tikzpicture}
			
		\end{minipage}
		\begin{minipage}{0.24\textwidth}
			\centering
			\begin{tikzpicture}
				[rotate around x=180]
				\foreach \i in {1,...,6}
				{
					\draw[rotate=-60*(\i)] (0, 1) node[circle, fill=white, draw=black!80, inner sep=0.7mm, minimum size=1mm]  (\i){};
					\draw (0, 0) node[circle, fill=white, draw=black!80, inner sep=0.7mm, minimum size=1mm]  (7){};
					\draw (0, -2) node[circle, fill=white, draw=black!80, inner sep=0.7mm, minimum size=1mm]  (8){};
				}
				
				\foreach \i/\j in {1/2, 2/3, 3/4, 4/5, 5/6, 6/1, 7/2, 7/4, 7/6, 8/3}
				{
					
					\draw[thick] (\i) -- (\j);
					\draw[bend right=60, thick] (8) to (1);
					\draw[bend left=60, thick] (8) to (5);
				}

			\end{tikzpicture}
			
		\end{minipage}
		\caption{Forbidden minors of the class of planar graphs whose duals are partial 3-trees.}
		\label{fig:forbiddenminorsofdualarepartial3trees}
		
	\end{figure}
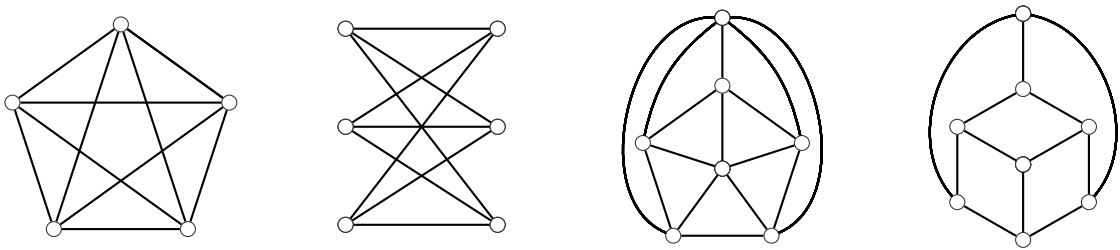

	A tight relation between this conjecture and restriction of \Cref{conj:NoOddK5ToSPC} to planar graphs is provided. We refer to \cite{NRS13} and references there in. In support of this conjecture we show that it holds for the subclass of planar multigraphs whose duals are partial 3-trees. In other words, graphs which do not contain any of the four graphs of \Cref{fig:forbiddenminorsofdualarepartial3trees} as minor. Before stating the claim we recall necessary tools from the literature.   
	
	The number of odd components in a graph $G$ is denoted by $o(G)$. Tutte's necessary and sufficient condition for a multigraph $G$ to have a perfect matching is as follows.
	
	\begin{theorem}\label{thm:Tutte's Theorem} A multigraph $G=(V,E)$ has a perfect matching if and only if $o(G-U)\leq |U|$ for any subset $U$ of vertices.
	\end{theorem}
	
	\begin{lemma}\label{lem:Perfectmatching} 
		Any $r$-regular graph $G$ with no odd cut of size smaller than $r$ has a perfect matching.		
	\end{lemma}
	
	\begin{proof}
		For any $U\subset V(G)$, let $G_1, G_2, \dots, G_l$ be the odd component(s) of $G-U$, where $l=o(G-U)$. On the one hand, there are at least $r$ edges connecting $G_i$ to $U$ for each $i\in\{1,2,\dots,l\}$. So the number of edges connecting $G_1\cup G_2\cup \dots \cup G_l$ to $U$ is at least $rl$. On the other hand, since $G$ is $r$-regular, there are at most $r|U|$ edges connecting $U$ to $G_1\cup G_2\cup \dots \cup G_l$. Therefore, $|U|\geq l=o(G-U)$. Hence one may apply Tutte's Theorem.		
	\end{proof}
	
	Given a planar graph together with an embedding on the plane, since facial cycles surrounded by a cycle $C$ determine the parity of $C$, we have: 
	
	\begin{observation}\label{obs:bipartite}
		If $G$ is $2k$-regular planar multigraph, then the dual of $G$ is bipartite.
	\end{observation}

	\begin{theorem} Any $r$-regular planar multigraph $G$ with no odd cut of size smaller than $r$ and whose dual is a partial 3-tree is $r$-edge colorable.
	\end{theorem}
	
	\begin{proof}
		The claim of the theorem for odd values of $r$ is proved in Theorem 7.4 of \cite{CN20}. Therefore, we only consider even values of $r$ and let $r=2k$.
		
		By \Cref{lem:Perfectmatching} $G$ has a perfect matching, let $M$ be a perfect matching of $G$. Let $H$ be the dual of $G$. We get that $H$ is bipartite by \Cref{obs:bipartite}. We assign the edges of $H$ corresponding to $M$ the negative sign, all other edges are assigned the positive sign. Denote the signed graph by $(H, \sigma)$. We claim that $(H, \sigma)$ has unbalanced girth at least $2k$.
		
		Since $G$ is $2k$-regular with a perfect matching $M$, we get that every face of $(H,\sigma)$ is a negative $2k$-cycle. Let $C$ be a negative cycle of $(H,\sigma)$ which is not a facial cycle. This cycle corresponds to an edge-cut $(X,X^c)$ of $G$. We claim that both $|X|$ and $|X^c|$ are odd. That is because vertices in $X$ correspond to the faces bounded by the cycle $C$. As the total number of negative faces of any signed plane graph is even, and as $C$ is a negative face in the subgraph induced by removing all elements not bounded by $C$, the number of negative faces that is bounded by $C$ is an odd number. Similarly, the number of faces not bounded by $C$ is an odd number. By the assumption on the number of edges of an odd edge-cut, which corresponds to the length of $C$, we conclude that $C$ is of length at least $2k$.     
		
		By the assumption, we know that $H$ is a bipartite partial 3-tree. Hence, $(H,\sigma)$ is a bipartite partial 3-tree of unbalanced girth $2k$. By \Cref{thm:Partial3Tree+Ok}, there is a homomorphism $\phi$ of $(H,\sigma)$ to $SPC(2k-1)$. In this mapping the image of every negative cycle contains a negative cycle. Since each face of $(H, \sigma)$ is a negative $2k$-cycle, and since the negative girth of $SPC(2k-1)$ is $2k$, the image of each face of $(H,\sigma)$ is a negative cycle of length $2k$. Recall that on each negative $2k$-cycle of $SPC(2k-1)$ each element of $S_{2k-1}$ appears exactly once. Thus if as color each edge of $G$ with an element of $S_{2k-1}$, the one to which its dual in $(H,\sigma)$ is mapped to by $\phi$, we have a proper $2k$-edge-coloring of $G$.  		
	\end{proof}

	\section{Algorithmic implications}\label{sec:algorithms}
	
	Our result provides an efficient algorithm for the following:
	Input is a signed graph $\widehat{B}$ and two integers $k$ and $t$. We may (efficiently) verify that $\widehat{B}$ induces no cycle from $\mathcal{O}_k$. Assuming it does not, the decision question is: 
	
	\begin{problem}
	Does every signed partial $t$-tree which induces no element of $\mathcal{O}_k$ map to $\widehat{B}$.
	\end{problem}
	
	To answer the question we first form a partial $\widehat{B}$-distance graph $(\widehat{B},w)$ whose edges are all the pairs each of which is contained in common negative $k$-cycle of $\widehat{B}$.
	
	Next we form the set $\mathcal{W}$ of all the $(t+1)$-cliques of $(\widehat{B}, w)$. Then we test if $\mathcal{W}\neq \emptyset$ and if it is $(t, k)$-closed. If $\mathcal{W}= \emptyset$, then we report that there exists a signed partial $t$-tree which induces no element of $\mathcal{O}_k$ and yet does not map to $\widehat{B}$. If $\mathcal{W}\neq \emptyset$ and it is, furthermore, $(t, k)$-closed, then we report that every signed partial $t$-tree which induces no element of $\mathcal{O}_k$  maps to $\widehat{B}$. If neither of the two cases happens, then we have found some cliques in $\mathcal{W}$ at which the condition of being $(t, k)$-closed fails. We modify $\mathcal{W}$ by removing all such cliques and repeat the process with the modified $\mathcal{W}$. Since $\mathcal{W}$ is finite, at a limited number of steps, either we reach a set of cliques that is $(t,k)$-closed, or a $\mathcal{W}$ which is empty.  
	When the conclusion is negative, that is we have reached $\mathcal{W}=\emptyset$, then we can trace back, and following the sequence by which the cliques are removed from $\mathcal{W}$, we can build an example of partial $t$-tree with no cycle from $\mathcal{O}_k$ which does not map to $\widehat{B}$.
	
	Noting that the study is motivated originally by potential extensions of the four-color theorem, the following algorithmic question is of high interest:
	
	\begin{problem}
		Given a positive integer $k$ and a signed graph $\widehat{B}$ which induces no element of $\mathcal{O}_k$ as a subgraph, is it true that every signed planar graph with no element of $\mathcal{O}_k$ as a subgraph maps to $\widehat{B}$? 
	\end{problem}
	
	With similarity of the question to some recent developments in refinement of Lov\'asz vector, it would not be a big surprise if the question turns out to be an undecidable one.
	
	Recall that Lov\'asz vector of a graph $G$ is a vector whose entries are labeled by all graphs (a countable number of them then) where the entry $H$ is the number of homomorphisms of $H$ to $G$. Lov\'asz showed that two (labeled) graphs are isomorphic if and only if they have the same Lov\'asz vector. 
	
	Restricting the vector by considering coordinates from a subclass of graphs, such as the subclass of trees, the subclass of partial $t$-trees, or the subclass of planar graphs, leads to various notions of similarities between graphs. We refer to \cite{Dev10}, \cite{MR20}, \cite{S24}, and references therein on the subject.
		
	It is shown in \cite{Dev10} that, for a fixed $t$, deciding if two graphs $G$ and $H$ admit the same number of homomorphisms from any partial $t$-tree can be done in polynomial time, and that an affirmative answer reveals certain similarities between the two graphs. For example, the case $t=1$ is about fractional isomorphism. 
	
	If instead we restrict the homomorphism count to the subclass of planar graphs, then it is shown in \cite{MR20} that two graphs have the same vector if and only if they are quantum isomorphic. However, in contrast to the previous example, it is shown that question of whether two given graphs are quantum isomorphic belongs to the family of undecidable problems.

	Our study can be viewed as the $0-1$ version of the Lov\'asz vector: an entry of the vector for $G$ is $1$ if the graph corresponding to that entry admits a homomorphism to $G$, $0$ otherwise. The density of homomorphism order (see \cite{HN04}), suggests that unless two graphs are homomorphically equivalent (i.e. have isomorphic cores), the $0-1$ version of Lov\'asz vector can distinguish them. In other words, if the $0-1$ version of Lov\'asz vector cannot separate the two graphs, then they are similar in the sense that they have the same core. The question then is what if the $0-1$ version of the vector is also restricted to the entries which are partial $t$-tree, or planar?
	
	 Restricting these $0-1$ vectors to coordinates corresponding to the class of partial $t$-trees, we classified signed graphs for which the associated vector is determined by the coordinates corresponding to the negative cycles.

	\section{Concluding remarks}\label{sec:remarks}
	
	For given integers $t$ and $k$, it is of high interest to find the smallest signed graph with no cycle from $\mathcal{O}_k$ which bounds the class $\mathcal{SPT}_{t,-k}$ (signed partial $t$-trees with no cycle from $\mathcal{O}_k$). 
	
	The example of bounds we have provided in this work for the case $t=3$, namely the signed projective cube of dimension $k-1$, are indeed optimal bounds. In \cite{NSS15} for each $k$ an example $G_k$ of a signed planar graph without cycles from $\mathcal{O}_k$ are given where each $G_k$ contains $2^{k-1}$ vertices identifying any pair of which would result in a cycle in $\mathcal{O}_k$. Upon further inspection of these examples, one may observe that they are, furthermore, partial 3-trees.
	
	Thus there is an implicit notion of prefectness for the class $\mathcal{SPT}_{3,-k}$: the minimum number of colors to color every member of $\mathcal{SPT}_{3,-k}$ in such a way that identifying vertices of the same color does not result in a cycle in $\mathcal{O}_k$ is the same as the minimum order of vertices of a signed graph $\widehat{B}$ with no cycle from $\mathcal{O}_k$ which bounds $\mathcal{SPT}_{3,-k}$. 
	
	This property does not hold for other values of $t$. In particular it is shown, in an unpublished work of W. He, the second author, and Q. Sun, that the maximum number of such set of vertices in a graph in $\mathcal{SPT}_{2,-3}$ is 6, whereas the minimum order of a homomorphism bound is 8. The smallest order of homomorphism bound for $\mathcal{SPT}_{2,-4}$ is 12 (unpublished work of the second author and R. Xu), for $\mathcal{SPT}_{2,-5}$ it is  15 \cite{BFN19} (see also \cite{BFN22}). For larger values of $k$, the problem remains open; perhaps the best order of a homomorphism bound for $\mathcal{SPT}_{2,-(2k+1)}$ is $(k+1)(k+3)$ and for $\mathcal{SPT}_{2,-2k}$ is $(k+1)(k+2)$. Bounds of quadratic order are given in \cite{FNX23}.       
		
	\section*{Acknowledgement} 
	
	This work has received support under the program ``Investissement d'Avenir" launched by the French Government and implemented by ANR, with the reference ``ANR‐18‐IdEx‐0001" as part of its program ``Emergence". Meirun Chen is supported by Fujian Provincial Department of Science and Technology (2024J011197).

\end{document}